\documentclass[11pt]{article}

\usepackage{graphicx}
\usepackage{amsmath}
\usepackage{amsfonts}
\usepackage{fancyhdr} 

\usepackage{wrapfig}

\newtheorem{lem}{Lemma}
\newtheorem{thm}{Theorem}

\newtheorem{prop}[lem]{Proposition}
\newtheorem{cor}[lem]{Corollary}

\newcommand{\sss}{\mathcal{S}}

\newcommand{\real}{\mathbb{R}}

\newcommand{\ddd}{\mathbb{D}}

\newcommand{\complex}{\mathbb{C}}

\newcommand{\aaa}{\mathcal{A}}

\newenvironment{pf}{\noindent {\em Proof}.\ \ }{\hspace*{\fill}\rule{.5ex}{1.4ex}\,}

\newcommand{\www}{\mathcal{W}}
\newcommand{\eee}{\mathcal{E}}

\newcommand{\kkk}{\mathcal{K}}
\newcommand{\ggg}{\mathfrak{G}}
\newcommand{\sep}{\mathfrak{V}}

\newcommand{\spn}{\mbox{span\,}}

\title{The structure of test functions that determine weighted composition operators}
\author{Peter~C.~Gibson and Mohammad~S.~Tavalla}

\begin{document}

\pagestyle{fancyplain}

\maketitle

\begin{abstract}
In the context of analytic functions on the open unit disk, a weighted composition operator is simply a composition operator followed by a multiplication operator.  The class of weighted composition operators has an important place in the theory of Banach spaces of analytic functions; for instance, it includes all isometries on $H^p$ $(p\neq 2)$.  Very recently it was shown that only weighted composition operators preserve the class of outer functions.  

The present paper considers a particular question motivated by applications: Which smallest possible sets of test functions can be used to identify an unknown weighted composition operator?  This stems from a practical problem in signal processing, where one seeks to identify an unknown minimum phase preserving operator on $L^2(\real_+)$ using test signals.

It is shown in the present paper that functions that determine weighted composition operators are directly linked to the classical normal family of schlicht functions.   The main result is that a pair of functions $\{f,g\}$ distinguishes between any two weighted composition operators if and only if there exists a zero-free function $h$ and a schlicht function $\sigma$ such that $\spn\{f,g\}=\spn\{h\sigma,h\}$.   This solves completely the underlying signal processing problem and brings to light an intriguing geometric object, the manifold of planes of the form $\spn\{h\sigma,h\}$.   As an application of the main result, it is proven that there exist compactly supported pairs in $L^2(\real_+)$ that can be used to identify minimum phase preserving operators.  
\end{abstract}

\section{Introduction}  
Let $\aaa$ denote the space of all analytic functions $f:\ddd\rightarrow\complex$ on the open unit disk,
endowed with the usual topology of uniform convergence on compact sets.  Given any two functions $\psi,\varphi\in\aaa$, where $\varphi(\ddd)\subset\ddd$, define the multiplication operator
$M_\psi:\aaa\rightarrow\aaa$ by the formula 
\[
M_\psi f=\psi f\quad(f\in\aaa),
\]
and the composition operator
$
C_\varphi:\aaa\rightarrow\aaa$ is by the formula 
\[
C_\varphi f=f\circ\varphi\quad(f\in\aaa).
\]
The present article is concerned with the class $\www$ of weighted composition operators, consisting of all operators of the form $M_\psi C_\varphi$.  The central problem at issue is to characterize the smallest possible sets of test functions in $\aaa$ that will identify an arbitrary operator $A\in\www$.   The problem may be stated precisely as follows.   

A set $X\subset\aaa$ {separates points} in a given family of operators $\www^\prime\subset\www$ if the mapping $A\mapsto A|_X$ is injective on $\www^\prime$.  Put another way, $X$ separates points in $\www^\prime$ if for every $A_1\neq A_2$ in $\www^\prime$, there exists a point $f\in X$ such that $A_1f\neq A_2f$.\\[5pt]
\noindent\fbox{
\parbox{.95\textwidth}{
\textbf{Main Problem}:\ \ Give a systematic representation of the set of all smallest possible sets $X\subset\aaa$ that separate points in $\www$. 
}
}\\[5pt]
The pair $\{1,z\}$ evidently separates points in $\www$, since given a non-zero operator $A=M_\psi C_\varphi\in\www$, the functions 
\[
\psi=A1\quad\mbox{ and }\quad\varphi=\frac{Az}{A1}
\]
completely determine $A$.  Moreover, one can see from this example that the smallest sets separating points in $\www$ consist of precisely two elements (see Proposition~\ref{prop-single}).   However, not every linearly independent pair works.  For instance, it can be shown that the pair $\{\cos z,\sin z\}$ fails to separate points in $\www$.  A priori, it is not clear why one linearly independent pair should work while another doesn't, and this is essentially what the main problem seeks to resolve.  As it turns out, the full family of separating pairs has a rather subtle structure; it is intimately connected with the classical family $\sss$ of schlicht functions that map the disk onto any conformally equivalent region in $\complex$, appropriately normalized.  The main result of the present paper is that a pair $\{f,g\}$ separates points in $\www$ if and only if there is a function $\sigma\in\sss$ and a function $h(z)=e^{zk(z)}$, where $k\in\aaa$, such that 
\[
\spn\{f,g\}=\spn\{h\sigma, h\}.
\]

\subsection{Overview of the paper}

The main problem originates from a practical question in signal processing; this connection is detailed in Section~\ref{sec-signal}.   

The theoretical importance of weighted composition operators in the context of complex functional analysis, including their emergence in some recent results, is briefly surveyed below, in Section~\ref{sec-background}.  

The solution to the main problem is worked out in Sections~\ref{sec-characterization} and \ref{sec-parameterization}.  The first step, carried out in Section~\ref{sec-characterization}, is to find a tractable characterization of pairs of functions that separate points in $\www$.  Such a characterization is given in Theorem~\ref{thm-EW}.   

The second step is to give an explicit parameterization---in other words, a constructive representation with no redundancy---for the family of pairs that satisfy the characterization.   This second step is carried out Section~\ref{sec-parameterization}, culminating in Theorem~\ref{thm-parameter}, which is the paper's main result.   

In Section~\ref{sec-application} the main result is applied to solve a non-trivial problem in signal processing.  The paper concludes in Section~\ref{sec-conclusion} with a brief indication of directions for further research.

\subsection{Background on weighted composition operators\label{sec-background}}

The books by Shapiro \cite{Sh:1993} and Cowan and MacCluer \cite{CoMacCl:1995} present an overview of weighted composition operators' connections to classical function theory and operator theory; Chapter~5 of Mart\'inez-Avenda\~no and Rosenthal's book \cite{MaRo:2007} gives a more recent survey of composition operators in the context of $H^2_\ddd$.  Hoffman's classic \cite{Ho:1962} serves as a useful reference for the theory of Hardy spaces, and in particular for the various characterization of outer functions.   For reference, a function $f\in H^1_\ddd$ is outer if it is not identically zero and
\[
\log|f(0)|=\frac{1}{2\pi}\int_0^{2\pi}\log|f(e^{i\theta})|\,d\theta.
\]
Equivalently, $f$ is outer if $f/||f||_1$ is an extreme point of the closed unit ball in $H^1_\ddd$.  

A fundamental role of weighted composition operators was discovered in 1964 by Forelli \cite{Fo:1964}, who proved that for all $p\neq 2$ in the range $1\leq p<\infty$, every isometry of $H^p_\ddd$ into itself is a weighted composition operator.  
The result does not hold for $p=2$, since in this case the Hilbert space structure gives rise to many additional isometries.  Forelli's paper builds on earlier work of De Leeuw, Rudin and  Wermer \cite{dLeRuWe:1960}; they showed that isometries of $H^1_\ddd$ \emph{onto} itself are weighted composition operators.   Their proof exploits the fact (see \cite{dLeRu:1958}) that extreme points of the unit ball of $H^1_\ddd$ are precisely the outer functions of norm 1, and so any surjective isometry must permute the outer functions.   

Very recently a related result that applies to all $H^p_\ddd$ spaces was established in \cite{GiLa:Polya2011}.  The idea is to abandon isometries and simply consider outer functions.  It follows immediately from \cite[Theorem~2]{GiLa:Polya2011} that any linear operator on $H^p_\ddd$ that maps the class of outer functions into itself is necessarily a weighted composition operator.  (Actually much more is true; any operator that maps polynomials having no zeros in $\ddd$ to functions that are zero free in a neighbourhood of a point is a weighted composition operator.  See \cite{GiLa:Polya2011} for details.)  The result is suggested by a theorem in \cite{GiLaMa:Outer2011}, wherein an operator is proved to be a weighted composition operator under the stronger hypothesis that all shifted outer functions in $H^2_\ddd$ be mapped into shifted outer functions.   (However the techniques used in \cite{GiLaMa:Outer2011} are inadequate to handle the case of \emph{un}shifted outer functions.)   

The main problem considered in the present paper stems from an application of Theorem~2 in \cite{GiLa:Polya2011} to minimum phase functions on $\real_+$; see \cite{GiLa:Id2011} for full details. 

\subsection{Connection to signal processing\label{sec-signal}}

The following signal processing considerations comprise the motivation behind the main problem.  

A function $f\in L^2(\real_+)$ is defined to be minimum phase if its Fourier-Laplace transform, which belongs to the Hardy space $H^2_{\complex_+}$, is outer.  
  Minimum phase functions are important in a range of applications, and in particular, in geophysical imaging.  This is partly because the (generally unknown) operator that models the transmission of an acoustic signal through a section of layered earth has a very special property: it maps minimum phase signals to delayed minimum phase signals.  Geophysical imaging entails determining this unknown operator using test signals.   The series of papers \cite{GiLaMa:Outer2011}, \cite{GiLa:Polya2011}, \cite{GiLa:Id2011} has culminated in a precise characterization of the linear operators on $L^2(\real_+)$ that preserve delayed minimum phase signals, thereby narrowing down the class of unknown operators that one seeks to identify.   It is shown in \cite{GiLa:Id2011} that an operator preserves delayed minimum phase only if it is conjugate via the Fourier-Laplace transform to a weighted composition operator on the right half plane $\complex_+$.   Furthermore, there is a particular pair of $L^2(\real_+)$ test functions that serves to identify any such operator,  namely the pair $\{f,f^\prime\}$, where $f(t)=te^{-t}$.   

This result raises a basic question that has direct implications for geophysical imaging: what other sets of test functions suffice to identify the class of operators in question?   This is especially important in analyzing data that has already been recorded, where the test signals are almost certainly different from those just mentioned.  It is known that not all pairs of linearly independent test functions can distinguish between minimum-phase preserving operators, and it is not a priori clear how to construct good test functions.   This question is easily translated to the setting of $H^2_\ddd$.  As will be seen in Section~\ref{sec-application}, solving the main problem in the broader setting of $\aaa\supset H^2_{\ddd}$, leads immediately to the desired Hardy space results, and with less technicalities than working directly in $H^2_{\ddd}$.  Thus the results in the present paper resolve completely the original signal processing question.

\section{Characterization of separating spaces\label{sec-characterization}}

The following simple observation streamlines the analysis considerably. 
\begin{prop}\label{prop-separation}
A set $X$ separates points in a family $\www^\prime$ if and only if the vector space $\spn X$ separates points in $\www^\prime$.  
\end{prop}
\begin{pf}  Given $A_1,A_2\in\www^\prime$, the statement $A_1f=A_2f$ holds for all $f\in X$ if and only if it holds for all $f\in\spn X$.  Therefore $X$ fails to separate points in a set $\www^\prime$ if and only if $\spn X$ fails to do so. 
\end{pf}

In what follows, it will be convenient to take advantage of Proposition~\ref{prop-separation} and to pass freely between sets and their spans when analyzing the separation of points in a family of operators.   Smallest possible separating sets are bases of minimum dimensional separating spaces.   Since all bases of a given space $V$ are in this sense equivalent, the more natural object to consider is $V$ itself, and not its bases.  Thus the main problem is more naturally stated as follows. \\[5pt]

\noindent\fbox{
\parbox{.96\textwidth}{
\textbf{Main Problem (reformulated)}:  Give a systematic representation of the set of all minimum dimensional subspaces $V\subset\aaa$ that separate points in $\www$. 
}
}\\[5pt]

Given a constant function $\varphi= z_0$, where $z_0\in\ddd$, the composition operator $C_\varphi$ will be denoted $C_{z_0}$.   This corresponds to evaluation at a point, since $C_{z_0}f=f(z_0)$, but where the value $f(z_0)\in\aaa$ has to be formally interpreted as a constant function and not a number.   (It will be clear from context whether an object is a number or a constant function, so this point will not be emphasized further.)  Let $\eee$ denote the subset of $\www$ consisting of all operators of the form 
$\alpha C_{z_0}$ where $\alpha\in\complex$ and $z_0\in\ddd$.   

(More abstractly, let $\kkk\subset\aaa$ denote the subspace comprised of constant functions.   Then $\eee$ may be characterized as the subset of operators in $\www$ mapping $\aaa$ into $\kkk$.  Since $\kkk\cong\complex$, this means that $\eee$ is essentially the intersection of the dual space of $\aaa$ with $\www$.)

Evidently a subspace $V\subset\aaa$ must separate points in $\eee$ if it is to separate points in $\www$; surprisingly, the converse is also true.  (This will be proved shortly---see Corollary~\ref{cor-EW}.)

\begin{prop}\label{prop-single}
No singleton $\{f\}\subset\aaa$ separates points in $\eee$. 
\end{prop}
\begin{pf}
The singleton $\{0\}$ clearly fails to separate points in $\eee$.  If $f$ is not identically 0, then there exist distinct points $z_0,z_1\in\ddd$ with $f(z_0)\neq0$.  Thus 
\[
\textstyle\frac{f(z_1)}{f(z_0)}C_{z_0}\neq C_{z_1},\quad\mbox{ but }\quad\textstyle\frac{f(z_1)}{f(z_0)}C_{z_0}f=C_{z_1}f.
\]
So $\{f\}$ does not separate points in $\eee$. \end{pf}

\begin{prop}\label{prop-E}
If a pair $\{f,g\}\subset\aaa$ separates points in $\eee$, then $f$ and $g$ have no common zero, and the meromorphic function $f/g:\ddd\rightarrow\complex\cup\{\infty\}$ is injective.   
\end{prop}
\begin{pf} Write $V=\spn\{f,g\}$, and let $z_0\in\ddd$.  If $f(z_0)=g(z_0)=0$, then $C_{z_0}|_V=0|_V$, but $C_{z_0}\neq 0$, so $V$ doesn't separate points in $\eee$.  Thus if $V$ separates points in $\eee$, then $f$ and $g$ have no common zero.  

Next suppose $f$ and $g$ have no common zero and consider $\mu=f/g$.   Suppose that $\mu(z_0)=\mu(z_1)$ for distinct points $z_0,z_1\in\ddd$.   Then $g(z_1)f(z_0)=f(z_1)g(z_0)$ and one of the four numbers $f(z_0),f(z_1),g(z_0),g(z_1)$ is non-zero.   Therefore either 
\[
g(z_1)C_{z_0}\neq g(z_0)C_{z_1}\quad\mbox{ or }\quad f(z_1)C_{z_0}\neq f(z_0)C_{z_1}.
\]
But both
\[
g(z_1)C_{z_0}|_V=g(z_0)C_{z_1}|_V\quad\mbox{ and }\quad f(z_1)C_{z_0}|_V=f(z_0)C_{z_1}|_V,
\]
so $V$ fails to separate points in $\eee$. Thus if $V$ separates points, $\mu$ is injective. 
\end{pf}\\
Note that by symmetry the conclusion of Proposition~\ref{prop-E} actually implies that both $f/g$ and $g/f$ are injective. 

\begin{thm}\label{thm-EW}
If $f,g\in\aaa$ have no common zero and the meromorphic function $f/g:\ddd\rightarrow\complex\cup\{\infty\}$ is injective, then $\{f,g\}$ separates points in $\www$. 
\end{thm}
\begin{pf}  As before, write $V=\spn\{f,g\}$ and $\mu=f/g$.  Suppose first that $A=M_\psi C_\varphi\in\www$ is arbitrary.  Then for all $z\in\ddd$ such that $\psi(z)\neq 0$, at least one of $Af(z)=\psi(z)f(\varphi(z))$ and $Ag(z)=\psi(z)g(\varphi(z))$ is different from zero, since $f$ and $g$ have no common zero.  So if $A\neq0$, then $A|_V\neq0$, showing that $V$ distinguishes any non-zero operator in $\www$ from the zero operator.  

Next let $A_1=M_{\psi_1}C_{\varphi_1}$ and $A_2=M_{\psi_2}C_{\varphi_2}$ be non-zero elements of $\www$, and suppose that $A_1|_V=A_2|_V$.  Since $A_1$ and $A_2$ are non-zero, the same is true of $\psi_1$ and $\psi_2$, so that the set $Z$ of points $z\in\ddd$ at which $\psi_1(z)\psi_2(z)=0$ is at most countable.  At each $z\in\ddd\setminus Z$, 
\[
\mu\circ\varphi_1(z)=\frac{A_1f(z)}{A_1g(z)}=\frac{A_2f(z)}{A_2g(z)}=\mu\circ\varphi_2(z).
\]
Since $\mu$ is injective it follows that $\varphi_1=\varphi_2$.   Since $f$ and $g$ have no common zeros, one of $f(\varphi_1(z))$ or $g(\varphi_1(z))$ is non-zero, for every $z\in\ddd\setminus Z$.  Thus at least one of $f\circ\varphi_1$ or $g\circ\varphi_1$ is non-zero on an uncountable set; suppose for definiteness that $f\circ\varphi_1$ has this property.   Then, given that $\varphi_1=\varphi_2$, 
\[
\psi_1(z)=\frac{(A_1f)(z)}{f\circ\varphi_1(z)}=\frac{(A_2f)(z)}{f\circ\varphi_2(z)}=\psi_2(z)
\]
for uncountably many $z\in\ddd$.  It follows by analyticity that $\psi_1=\psi_2$ and hence that $A_1=A_2$, proving that $V$ separates points in $\www$.   \end{pf}\newline
In conjunction with Proposition~\ref{prop-E}, Theorem~\ref{thm-EW} immediately yields the following corollary. 
\begin{cor}\label{cor-EW} A pair $\{f,g\}\subset\aaa$ separates points in $\www$ if and only if it separates points in $\eee$.   
\end{cor}

\section{Parameterization of separating spaces\label{sec-parameterization}}

Given a zero-free function $g\in\aaa$ and a schlicht function $\sigma$, the pair $\{g\sigma,g\}$ satisfies the hypothesis of Theorem~\ref{thm-EW}, and therefore the vector space $V=\spn\{g\sigma,g\}$ separates points in $\www$.   Moreover, it is minimum-dimensional.   
The next theorem asserts that every two-dimensional space $V$ separating points in $\www$ arises this way.   To represent them in a systematic way, a bijective parameterization of the family of all minimum dimension (i.e. two-dimensional) separating spaces will be constructed.  This requires some setup, as follows. 

Let $\ggg_2$ denote the family of all two-dimensional subspaces of $\aaa$, and let $\sep_2\subset\ggg_2$ denote the subfamily consisting of spaces $V$ that separate points in $\www$.  (Roughly speaking, $\ggg_2$ is like an infinite-dimensional grassmannian, with $\sep_2$ being a special submanifold to be parameterized.  Indeed one can formalize this perspective.  However, putting charts on $\ggg_2$ requires cutting down from $\aaa$ to the Hardy space $H^2_{\ddd}$, and then using the Hilbert space structure of $H^2_{\ddd}$.  The geometric perspective will not be pursued in the present article---although it provides a useful picture.)  

Let $\aaa_\ast\subset\aaa$ denote the multiplicative group consisting of all $f\in\aaa$ such that: (i) $0\not\in f(\ddd)$; and (ii) $f(0)=1$.  Let $\sss\subset\aaa$ denote the classical family of schlicht functions, that is, the family of all injective analytic mappings $\sigma:\ddd\rightarrow\complex$ with the normalization $\sigma(0)=0$ and $\sigma^\prime(0)=1$.   
For present purposes a further reduction of $\sss$ is needed.   Define $\sigma,\tau\in\sss$ to be equivalent, writing $\sigma\sim\tau$, if they are related by an automorphism of the Riemann sphere.  To be explicit, $\sigma\sim\tau$ if there exists a point $\xi\in\complex\cup\{\infty\}$ such that 
\[
\tau=\frac{\sigma}{1-\frac{1}{\xi}\sigma}\quad\mbox{ equivalently }\quad \sigma=\frac{\tau}{1+\frac{1}{\xi}\tau}.
\]
Note that since $\tau\in\sss$, it is required that $\xi\not\in\sigma(\ddd)$; on the other hand, $\xi=\infty$ is permitted, whence $\sigma\sim\sigma$.  (The relation $\sim$ is easily seen to be transitive and hence an equivalence relation; furthermore, equivalence classes are closed.)  Let $\widetilde{\sss}$ denote a section of the quotient structure $\sss/\sim$ (which inherits compactness from $\sss$), meaning that $\widetilde{\sss}\subset\sss$ contains precisely one representative of each $\sim$ equivalence class in $\sss$.    

\begin{thm}\label{thm-parameter}
The mapping $\Phi:\aaa_\ast\times\widetilde{\sss}\rightarrow\ggg_2$ defined by the formula
\[
\Phi(f,\sigma)=\spn\{f\sigma,f\}
\]
is a bijection from $\aaa_\ast\times\widetilde{\sss}$ onto $\sep_2$.   
\end{thm}
\begin{pf}
If $(g,\sigma)\in\aaa_\ast\times\widetilde{\sss}$, then the pair $\{g\sigma,g\}$ satisfies the hypothesis of Theorem~\ref{thm-EW} with $f=g\sigma$, and so $V=\spn\{g\sigma,g\}\in\sep_2$.   Thus 
\[
\Phi(\aaa_\ast\times\widetilde{\sss})\subset\sep_2.
\]
Next write $V_1=\Phi(g_1,\sigma_1)$ and  $V_2=\Phi(g_2,\sigma_2)$, and suppose that $V_1=V_2$.  Observe that, by the normalizations that define $\aaa_\ast$ and $\sss$, each of $V_1$ and $V_2$ contains a unique element, $g_1\sigma_1$ and $g_2\sigma_2$ respectively, which takes the value 0 at 0 and has derivative 1 at 0.   Therefore 
\begin{equation}\label{gsigma}
g_1\sigma_1=g_2\sigma_2.
\end{equation}
Since $g_1\in V_2$, and $g_1(0)=1$, it follows that $g_1$ has the form 
\[
g_1=g_2+\alpha g_2\sigma_2=g_2+\alpha g_1\sigma_1,
\]
and hence $g_2=g_1(1-\alpha\sigma_1)$.   It follows in turn from (\ref{gsigma}) that 
\[
\sigma_2=\frac{\sigma_1}{1-\alpha\sigma_1},
\]
so that $\sigma_2\sim\sigma_1$.   By definition of $\widetilde{\sss}$, this forces $\sigma_2=\sigma_1$, which in turn by (\ref{gsigma}) forces $g_1=g_2$, proving that $\Phi$ is injective.  It remains only to show that every point of $\sep_2$ belongs to the range of $\Phi$.  

Let $V=\spn\{f,g\}\in\sep_2$ be an arbitrary two-dimensional space that separates points in $\www$.  By Proposition~\ref{prop-E}, $f$ and $g$ have no common zero, and the meromorphic functions $f/g$ and $g/f$ are injective.  One of $f(0)$ and $g(0)$ is different from zero; suppose for definiteness that $g(0)\neq0$, and set $\mu=f/g$.  Let 
\[
\alpha\in\complex\setminus\mu(\ddd).
\]
Such an $\alpha$ is guaranteed to exist, since the disk is not conformally equivalent to the Riemann sphere, or the Riemann sphere with a point deleted.  Set 
\[
\sigma=\lambda\frac{\mu-\mu(0)}{\mu-\alpha}=\lambda\frac{f-\mu(0)g}{f-\alpha g},
\]
where $\lambda\in\complex$ is chosen so that $\sigma^\prime(0)=1$.   Then $\sigma\in\sss$, and the function $G=(f-\alpha g)/(f(0)-\alpha g(0))$ belongs to $\aaa_\ast$.   Furthermore,
\[
\spn\{G\sigma,G\}=\spn\{f,g\}.
\]
Finally, let $\tau\in\widetilde{\sss}$ be equivalent to $\sigma$, with 
\[
\tau=\frac{\sigma}{1-\frac{1}{\xi}\sigma}.
\]
Then, setting $h=(1-\frac{1}{\xi}\sigma)G$, it follows that $(h,\tau)\in\aaa_\ast\times\widetilde{\sss}$, and 
\[
\spn\{h\tau,h\}=\spn\{G\sigma,G\}=\spn\{f,g\}=V,
\]
proving that $V$ is in the range of $\Phi$.  \end{pf}\\

Note that the group $\aaa_\ast$ consists of all functions of the form $f(z)=e^{zg(z)}$, where $g\in\aaa$ is arbitrary.  Thus is it easy to produce such functions.   On the other hand, while univalent functions have a long history (see \cite{Du:1983}), explicitly constructing schlicht functions $\sigma\in\sss$ is not something that can be done systematically.  De Branges's Theorem \cite{dBr:1985} gives necessary conditions in terms of the Taylor coefficients but there is no known necessary and sufficient condition. 

\section{An application to signal processing\label{sec-application}}

The signal processing perspective outlined in Section~\ref{sec-signal}  was the original impetus for the main problem, and it is a source of numerous particular questions.   For example, do there exist compactly supported test functions in $L^2(\real_+)$ that can be used to identify an unknown minimum phase preserving operator on $L^2(\real_+)$?   (Any physically realizable test function must be compactly supported.) 

This question can be solved in the affirmative with the help of Theorem~\ref{thm-parameter}, as follows.   Note that one may pass between the Hardy spaces on the half plane and the disk by way of the outer preserving isometry 
\[
\Phi:H^2_{\complex_+}\rightarrow H^2_{\ddd},\quad (\Phi f)(z)=\frac{2\sqrt{\pi}}{(1+z)}f\Bigl(\frac{1-z}{1+z}\Bigr).
\]
(See \cite{Ho:1962} and \cite{GiLa:Id2011}.)  Let $\mathcal{F}$ denote the Fourier-Laplace transform
\[
\mathcal{F}:L^2(\real_+)\rightarrow H^2_{\complex_+},\quad (\mathcal{F}g)(z)=\frac{1}{\sqrt{2\pi}}\int_{0}^\infty g(t)e^{-zt}\,dt,
\]
which is also an isometry.  An operator $A$ on $L^2(\real_+)$ preserves minimum phase if and only if 
\[
B=\Phi\mathcal{F}A(\Phi\mathcal{F})^{-1}:H^2_\ddd\rightarrow H^2_\ddd
\]
maps outer functions to outer functions, which by \cite[Theorem~2]{GiLa:Polya2011} is only possible if $B$ is a weighted composition operator.  A pair $\{f,g\}$ of test functions determines $A$ if and only if their image $\{ \Phi\mathcal{F}f,\Phi\mathcal{F}g\}$ determines $B$.  Note that the family of operators $\eee$ is comprised of bounded operators on $H^2_\ddd$ (technically each operator in $\eee$ restricts to a bounded operator on Hardy space), and except for 0, each such operator preserves outer functions (since non-zero constant functions are outer).   Therefore, if a pair of test functions is to identify $B$, it must, according to Corollary~\ref{cor-EW}, separate points in $\www$.   

Hence to construct a pair $\{f,g\}$ of compactly supported test functions that determines $A$ it suffices that their image $\{ \Phi\mathcal{F}f,\Phi\mathcal{F}g\}$ be generated according to Theorem~\ref{thm-parameter}.   Consider first the compactly supported function
\begin{equation}\label{g}
g(t)=\chi_{[0,2\pi]}(t)e^{-t}\sin t,
\end{equation}
which has Fourier-Laplace transform
\[
(\mathcal{F}g)(z)=\frac{1-e^{-2\pi(1+z)}}{\sqrt{2\pi}(1+(1+z)^2)}.
\]
The function $\mathcal{F}g$ is zero-free on $\complex_+$, and it follows that $\rho=\Phi\mathcal{F}g\in H^2_\ddd$ is zero-free on $\ddd$.   Thus by Theorem~\ref{thm-parameter}, the $H^2_\ddd$ pair $\{\rho\sigma,\rho\}$ separates points in $\www$ for any schlicht function $\sigma$.   The simplest choice is to take the identity $\sigma(z)=z$, then define
\[
f=\mathcal{F}^{-1}\Phi^{-1}(\rho\sigma),
\]
and check that $f$ has compact support.   Note that $\Phi$ is self-inverse, up to a constant:
\[
\Phi^{-1}:H^2_{\ddd}\rightarrow H^2_{\complex_+},\quad (\Phi k)(z)=\frac{1}{\sqrt{\pi}(1+z)}k\Bigl(\frac{1-z}{1+z}\Bigr).
\]
Thus 
\[
(\mathcal{F}f)(z)=\Phi^{-1}(\rho\sigma)(z)=(\mathcal{F}g)(z)\cdot\frac{1-z}{1+z}.
\]
The inverse Fourier-Laplace transform of $(1-z)/(1+z)$ is the tempered distribution 
\[
\nu(t)=\sqrt{2\pi}\bigl(2\chi_{[0,\infty)}(t)e^{-t}-\delta(t)\bigr).
\]
Therefore, 
\begin{equation}\label{f}
\begin{split}
f(t)&=\frac{1}{\sqrt{2\pi}}(g\ast \nu)(t)\\
&= \left(2\int_0^{\min\{t,2\pi\}}e^{-(t-x)}e^{-x}\sin x\,dx\right)-g(t)\\
&=\left(2e^{-t}\int_0^{\min\{t,2\pi\}}\sin x\,dx\right)-g(t)\\
&=\chi_{[0,2\pi]}(t)e^{-t}(2-2\cos t-\sin t).
\end{split}
\end{equation}
Thus the pair of test functions $\{f,g\}$ defined by (\ref{f}) and (\ref{g}) serves to identify minimum phase preserving operators on $L^2(\real_+)$, proving that it is possible for such functions to have compact support.  (The pair $g,\nu$ is rather exotic in that its convolution has the same support as one of the factors; this is the reason for choosing $g$ as a product of an exponential and a function that integrates to zero over the support.)  Of course one may replace $\{f,g\}$ by any pair having the same span.   For instance, setting $h(t)=\chi_{[0,2\pi]}(t)e^{-t}(1-\cos t)$, one sees that 
\[
\spn\{h,h^\prime\}=\spn\{f,g\},
\]
so the compactly supported pair $\{h,h^\prime\}$ also serves to identify minimum phase preserving operators.  As a concluding remark, note that the occurrence of the derivative is no accident.   Differentiation in $L^2(\real_+)$ corresponds to multiplication by $z$ in $H^2_{\complex_+}$; since the function $z$ is univalent on $\complex_+$, this engenders a pair that separates points in the family of weighted composition operators as long as the Fourier-Laplace transform of the original function is zero-free on $\complex_+$.    

\section{Conclusions\label{sec-conclusion}}

While the identification of weighted composition operators using test functions comes from a practical question in signal processing, the parametrization of all two-dimensional separating spaces that emerges (Theorem~\ref{thm-parameter}) is of considerable mathematical interest in its own right.  By virtue of Corollary~\ref{cor-EW} and the fact that scaled point evaluation is a continuous operator on $H^p_\ddd$ for every $1\leq p\leq\infty$,  Theorem~\ref{thm-parameter} is easily adapted to any of the Hardy spaces on the disk.  The corresponding sets
\[
\sep_2^p=\{V\in\sep_2\,|\,V\subset H^p_{\ddd}\}
\]
comprise an intriguing family whose geometric structure calls for further investigation.  The most interesting case is perhaps $\sep_2^2$, where one can exploit the Hilbert space structure to mimic the usual way of defining charts for finite-dimensional grassmannians.  Given that such an infinite-dimensional  construction arises in the theory of classifying spaces, it is of interest to know what $\sep_2^2$ looks like from the topological standpoint---see \cite{Le:2011}.

There is a great deal more to be said concerning the application of Theorem~\ref{thm-parameter} to signal processing; the example discussed in Section~\ref{sec-application} is not meant to be exhaustive.   Rather, it is intended merely to illustrate the utility in a practical setting of the paper's main result.  Further analysis  of minimum phase preserving operators will be dealt with in a separate paper.   

\bibliography{ReferencesWC}

\end{document}